\documentclass{amsart}

\theoremstyle{definition}
\newtheorem{ex}{Example}
\newtheorem{rem}{Remark}

\theoremstyle{plain}
\newtheorem{defin}{Definition}
\newtheorem{thm}[defin]{Theorem}
\newtheorem{lem}[defin]{Lemma}

\newcommand{\CN}{\mathbb{C}^{N}}
\newcommand{\CNp}{\mathbb{C}^{N^\prime}}
\newcommand{\Cn}{\mathbb{C}^{n}}
\newcommand{\Cnp}{\mathbb{C}^{n^\prime}}
\newcommand{\Cd}{\mathbb{C}^{d}}

\newcommand{\Rd}{\mathbb{R}^{d}}

\newcommand{\C}{\mathbb{C}}
\newcommand{\crb}[1]{\mathcal{V} (#1)}
\newcommand{\Nn}{\mathbb{N}^{n}}

\newcommand{\vardop}[2]{\frac{\partial#1}{\partial
#2}}

\DeclareMathOperator{\spanc}{span _\C}
\DeclareMathOperator{\imag}{Im}
\DeclareMathOperator{\real}{Re}
\DeclareMathOperator{\codim}{codim}

\begin{document}
\title[A reflection principle]{A reflection principle
for  real--analytic submanifolds of complex spaces}
\date{April 20, 1999}
\author{Bernhard Lamel}
\address{Department of Mathematics\\
University of California San Diego\\
La Jolla, CA 92092--0112}
\email{blamel@euclid.ucsd.edu}
\begin{abstract}
We give an invariant nondegeneracy condition
for CR--maps between  generic submanifolds in
different dimensions and use it to prove a reflection
principle for these maps.
\end{abstract}
\maketitle
\section{An invariant condition for CR--maps}

Assume that $M\subset\CN$ and $M^\prime\subset \CNp$
are generic real--analytic submanifolds of $\CN$ and
$\CNp$, respectively. We will consider CR--maps
$H:M\to M^\prime$, that is, maps which push complex
tangent vectors to $M$ to complex tangent vectors to
$M^\prime$, ie, $H_{*}T^c_p M\subset T^c_{H(p)}
M^\prime$. For our purposes, this means that after
choosing coordinates $(Z_1^\prime,\dots,
Z_{N^\prime}^\prime)\in\CNp$ the components to $H$
in this coordinate system are CR--functions on $M$
(For details concerning CR--maps and the last
statement, the reader is encouraged to see
\cite{BERbook}, \S 2.3.). We want to give a
condition on the map $H$ which guarantees that it
is the restriction of a holomorphic map from $\CN$
into $\CNp$. Our nondegeneracy condition involves
taking derivatives of the complex gradients to
$M^\prime$ pulled back to $M$ via $H$ with respect
to the CR--vector fields tangent to $M$. For
simplicity, assume that $0\in M$, $0\in M^\prime$
and $H(0)= 0$ (in the examples, this will not always be the case). Let
$\rho^\prime =
(\rho^\prime_1,\dots,\rho^\prime_{d^\prime})$ be a
local defining function for $M^\prime$ near
$0\in\CNp$. We write $\crb{M}$ for the
CR--bundle of $M$.
\begin{defin} \label{D:nondeg} Assume that $H$ is a
CR--mapping of $M$ into $M^\prime$ of class
$C^{k_0}$ and $H(0)=0$. Define complex linear
subspaces
$E_k (0)\subset\CNp$ by
\begin{multline} \label{E:ek}
E_k(0) = \spanc \{ L_1
\cdots L_j \,
\rho^\prime_{s Z^\prime}
 \left(H(Z),\overline{H(Z)}\right)\Bigr|_{Z=0}
\colon L_r
\in
\Gamma( M ,\crb{M}) \\ \text{ for } 1\leq r \leq
j , 0\leq j \leq k, 1\leq l\leq d^\prime
\}.
\end{multline}
We say that $H$ is $k_0$--nondegenerate if $E_{k_0}(0)
= \CNp$ and $E_k(0) \neq \CNp$ for $k < k_0$. Here we
have used the notation $\rho^\prime_{s Z^\prime} = (
\rho^\prime_{s Z^\prime_1},\dots,
\rho^\prime_{s Z^\prime_{N^\prime}}
)$ for the complex gradient.
\end{defin}
 We
will show that Definition~\ref{D:nondeg} is invariant
under biholomorphic changes of coordinates in Lemma
\ref{L:ind}.
\begin{rem} If there exists a $k_0$--nondegenerate
map into the manifold $M^\prime$, then automatically
$M^\prime$ is finitely nondegenerate. Here we say
that a real submanifold $M$ is finitely
nondegenerate (or more specifically
$k_0$--nondegenerate) if the identity mapping of
$M$ is $k_0$--nondegenerate in the sense of
Definition~\ref{D:nondeg}. For the original definition
of finite nondegeneracy, see
\cite{BERbook}, \S 11.1.. In fact, if there
is a $k_0$--nondegenerate map $H$ into
$M^\prime$, then $M^\prime$ will be
$\ell_0$--nondegenerate with
$k_0\geq\ell_0$. The reader can check this by
writing out the condition in Definition
\ref{D:nondeg} in coordinates and using the chain
rule.
\end{rem}
Let $\Gamma$ be an open, convex cone in $\Rd$.
By a wedge $\hat W_\Gamma\subset\CN$with edge $M$ we
mean a set of the form $\hat W_\Gamma=\{Z\in U \colon
\rho(Z,\bar Z )
\in
\Gamma \}$ where $\rho$ is a defining function for
$M$ and
$U\subset\CN$ is an open neighbourhood of the origin.
We can now state our result.
\begin{thm} \label{T:main}
Let $M\subset\CN$, $M^\prime\subset\CNp$ be generic
real--analytic submanifolds, $0\in M$, $0\in
M^\prime$. Assume that
$H:M\to M^\prime$ is a CR--map with $H(0)=0$ which is
$k_0$--nondegenerate at $0$ and extends
continuously to a holomorphic function (called
again $H$) in a wedge
$\hat W_\Gamma$ with edge
$M$. Then $H$ extends to a holomorphic map
$H:\CN\to\CNp$ in a full neighbourhood of $0\in\CN$.
\end{thm}

Let us note that if we assume $M$ to be of finite
type, then by a result of Tumanov \cite{TU1} every
CR--map
$H$ will automatically extend to a wedge. Also, the
hypotheses that $H$ extends continuously and that $H$
is
$C^{k_0}$ on $M$ implies that actually $H$ extends as
a
$C^{k_0}$ function up to the edge $M$
of $\hat W_\Gamma$ (see e.g. \cite{BERbook},
\S 7.5.).

Note that if
$M\subset\CN$,
$M^\prime\subset\CN$ are of equal codimension and both
$\ell_0$--nondegenerate, then every
CR--diffeomorphism of class $\ell_0$ is
$\ell_0$--nondegenerate. Let us now illustrate
Definition \ref{D:nondeg} by giving some examples.

\begin{ex}\label{EX:first}
Let $M=\{(z,w)\in\C^2 \colon \imag w = |z|^4 \} $,
$M^\prime = \{(z,w)\in\C^2 \colon \imag w = |z|^2 \}$.
Then the map $H(z,w)=(z^2,w)$ is $2$--nondegenerate.
\end{ex}
We will now give an example where $M$ has
codimension $1$ and $M^\prime$ has codimension
$2$.
\begin{ex}
Consider the manifolds
$M=\{(z,w)\in\C^2
\colon
\imag w = |z|^2 \}$,
$M^\prime = \{ (z^\prime,w_1^\prime,w_2^\prime)\colon
\imag w_1^\prime = \imag w_2^\prime = |z^\prime |^2
\}$. The map $H:\C^2\to\C^3$, $H(z,w) =
(z,w,w)$ is $1$--nondegenerate.\end{ex}
\begin{ex}
As the source manifold consider the
ball in $\C^2$, $M= \{ (z_1, z_2)\in \C^2 \colon
|z_1|^2 + |z_2|^2 = 1 \}$. There are (up to
automorphisms) exactly three nonlinear maps from $M$
into the $3$--ball
$M^\prime = \{ (z_1^\prime, z_2^\prime,
z_3^\prime)\in \C^3
\colon  |z_1^\prime|^2+ |z_2^\prime|^2 +
|z_3^\prime|^2 = 1
\}$, as Faran \cite{FA1} showed: $H_1(z_1,z_2) =
(z_1,z_1 z_2, z_2^2)$,  $H_2(z_1,z_2) =
(z_1^2,\sqrt{2}z_1 z_2, z_2^2)$, and $H_3(z_1,z_2) =
(z_1^3,\sqrt{3}z_1 z_2, z_2^3)$. These are $2$,
$2$, and $3$--nondegenerate, respectively, at the
point
$(1,0,0)$. However, if we consider maps from an
$n$--ball into an $n+1$--ball which are at least
$C^3$ for
$n\geq 3$, a result due to Webster \cite{WE1} tells
us that these are all linear and hence never
nondegenerate.
\end{ex}

We would like to note that variants of
Definition \ref{D:nondeg} have appeared in the
literature before, for example in the case of
pseudoconvex hypersurfaces (\cite{CSK}, \cite{FA2})
but the invariant
formulation given here is new.  Another similar
condition can be found in \cite{HAN}. In this paper
Han studies the reflection principle for
Levi--nondegenerate hypersurfaces. He requires that
after putting the target hypersurface in suitable
coordinates $(z^\prime,w^\prime)$ and writing
$H=(f,g)$ in these coordinates, that the
derivatives $L^\alpha \bar f (0)$ where $L^\alpha =
L_1^{\alpha_1}\cdots L_n^{\alpha_n}$ and
$L_1,\dots,L_n$ is a basis for the CR--vector
fields tangent to $M$ span $\Cnp$ (where $N^\prime=
n^\prime+1$). The problem with this condition is
that it is not invariant under biholomorphic
changes of coordinates, and Theorem 1.1 in
\cite{HAN} is incorrect, as the following
example which is due to D. Zaitsev shows.

\begin{ex}\label{EX:dime}
Consider the hypersurface $M^\prime \subset \C^3$
given by
$\imag w = |z_1|^2 - |z_2|^2$. After the change of
coordinates $z_1 = \zeta_1 +\zeta_2 - \zeta_2^2$,
$z_2 = \zeta_2$, $w = \tau$, $M$ takes the form
$\imag \tau = |\zeta_1 +\zeta_2 -\zeta_2^2|^2 -
|\zeta_2 |^2$. Let $\phi$ be a CR--function on $M=
\{ (z,w) \colon \imag w = |z|^2 \}$ and  consider the
map $H(z,w)= (\phi^2 (z,w), \phi (z,w), 0)$.Then $f
= (\phi^2,\phi)$ and if we choose $\phi$ such that
$\phi_z(0)\neq 0$, then the derivatives $L \bar f(0)$
and
$L^2 \bar f (0)$ will span $\C^2$. Choosing for
$\phi$ any CR--function which does not extend
(such exist since
$M$ is strongly pseudoconvex) gives a counterexample
to Theorem 1.1 in \cite{HAN}.
\end{ex}

We would now like to mention some of the reflection
principles which have been obtained for manifolds in
complex spaces of different dimensions. In the case
of pseudoconvex hypersurfaces,  reflection principles
 were obtained in \cite{FO1} and
\cite{HU1}. A reflection principle for balls has
been obtained in \cite{BHR}.
Also, the reader is referred to the survey article
\cite{FO2} for some other special cases which have
been settled. 

 We will conclude this
section by proving that Definition~\ref{D:nondeg} is
independent of the choice of coordinates. The
definition is independent of the choice of defining
function, as the reader can easily prove by induction
on $k$. The definition is also independent of the
choice of coordinates
$Z\in\CN$. It only remains to be shown that it is
independent of  changes of holomorphic
coordinates in the target space.
\begin{lem}\label{L:ind}
Let $M\subset\CN$, $M^\prime\subset\CNp$, $H$ be as in
Definition \ref{D:nondeg}. Define $E_k(0)$ by
\eqref{E:ek}. If we change coordinates by $\tilde
Z^\prime = F(Z^\prime)$ in $\CNp$, and denote the
corresponding spaces defined by \eqref{E:ek} in the
new coordinates by $\tilde E_k(0)$, then
\begin{equation}\label{E:ektransform}
\tilde E_k (0) = E_k (0) \left(
\vardop{F}{Z^\prime}(0)\right)^{-1},
\end{equation}
where  $E_k(0)$ is considered as a space of row
vectors.
\end{lem}
\begin{proof}
Suppose that
$F:\CNp\to\CNp$ is a local biholomorphic change of
coordinates with
$F(0)=0$. Since Definition~\ref{D:nondeg} is
independent of the choice of defining function we
can choose
$\tilde\rho =
\rho^\prime
\circ F^{-1}$ as a defining function
 for
$F(M)$.  By the chain rule we have that
\[\tilde \rho_{ j \tilde Z^\prime}= (\rho^\prime_j
\circ F^{-1})_{\tilde Z^\prime} = (\rho^\prime_{j
Z^\prime}
\circ F^{-1})
\left( \vardop{F^{-1}}{\tilde{Z}^\prime} \right).
\]
Now the only crucial point is that the matrix on
the right hand side has holomorphic entries.
Hence, if we compose with $H$, this matrix will
be annihilated by the CR--vector fields on $M$. Applying
the vector fields
$ L_1,\dots, L_r$ as in \eqref{E:ek} to this
equation we conclude that
\begin{multline} L_1\cdots
L_r \tilde \rho_{ j \tilde Z^\prime} (F\circ
H(Z),\overline{F\circ H(Z)})
\\=  L_1
\cdots L_r
\rho^\prime_{j Z^\prime} (H(Z),\overline{H(Z)})
\left( \vardop{F^{-1}}{\tilde{Z}}(F(H(Z))) \right).
\end{multline}
Evaluating this identity at $0$, we obtain
\eqref{E:ektransform}.
\end{proof}
\section{Preliminaries}\label{S:prel}
Before we give the proof of
Theorem~\ref{T:main}, we want to give the
technical details on which the proof is based.
We will be making use of the Edge--of--the--Wedge
Theorem. The technique used here is basically
the same as in \cite{BJT}. We start by fixing a
defining equation for $M$. Assume that $\codim
M = d$. Then it is possible to choose
coordinates $(z,w)\in \Cn\times\Cd = \CN$ near
$0$ and a real--analytic function
$\varphi:\Cn\times\Rd \to \Rd$, defined near $0$, such
that
$M$ is given by
\begin{equation}\label{E:definingequM}
\imag w = \varphi (z,\bar z, \real w ),
\end{equation}
and $\varphi(0)=0$, $d\varphi (0)=0$, $\varphi
(z, 0,s) =
\varphi (0,\bar z,s) = 0$. Define the function
$\Psi : \Cn\times\Rd\times\Rd \to \CN$,
\[
\Psi(z,\bar z,s,t) = (z,s+it + i\varphi (z,\bar
z, s+ it) ).
\]
This is a real--analytic diffeomorphism, but
{\em not} holomorphic. It flattens out $M$ in the
sense that near the origin $\Psi^{-1} (M) =
\Cn\times\Rd\times
\{0\}$. With respect to $s+it$, that is, for $z$
fixed, $\Psi$ is actually holomorphic.
We shall equip $\Cn\times\Rd$
with the CR--structure transported by
$\Psi$ from $M$. Let
$\Gamma$ be an open, convex cone in
$\Rd$ (all cones will be assumed to be open
and convex). By a wedge
$W_\Gamma$ we  mean a set of the
form
\[
W_\Gamma = \{ (z,s,t) \in \Cn\times
\Rd\times\Rd \colon t\in\Gamma \}.
\]
We will also call the image under $\Psi$ of
$W_\Gamma$ a wedge, but to avoid confusion we
will write $\tilde W_\Gamma$ for it. We will
also write $W_\Gamma^\epsilon = \{ (z,s,t)\in
W_\Gamma \colon |z|<\epsilon, |s|<\epsilon, |t|
< \epsilon \}$ and $\tilde W_\Gamma^\epsilon$ for
the image under $\Psi$ of $W_\Gamma^\epsilon$. The
apparent difference with the definition of the
wedge $\hat W_\Gamma$ given before the statement of
Theorem~\ref{T:main} is reconciled by the
observation that given any cone
$\Gamma^\prime$ whose closure is contained in
$\Gamma$, then in the sense of germs of sets, $\hat
W_{\Gamma^\prime} \subset \tilde W_\Gamma$ and
$ \tilde W_{\Gamma^\prime} \subset \hat
W_\Gamma $ (see the end of Section 2 of
\cite{BJT}).

By saying that a CR--function
$h$ on $M$ extends holomorphically to
$W_\Gamma^\epsilon$ we shall mean that there
is a holomorphic function $f$ on
$\tilde W_\Gamma^\epsilon$, continuous up to $M$, such
that
$h(z,\bar z, s) = f(\Psi(z,\bar z,s,0))$, and
likewise we will say that $h$ extends
holomorphically to a full neighbourhood of the
origin if there is a holomorphic function $f$
defined in a full neighbourhood of the origin such
that this equality holds. To say that
$h$ extends holomorphically to a full
neighbourhood of the origin is the same as
having an
$\epsilon>0$ such that for each $z$ with
$|z|<\epsilon$ the function
$s\mapsto h(z,\bar z,s)$ extends
holomorphically in $s+it$ for $|s|<\epsilon$,
$|t|<\epsilon$ such that the function $h(z,\bar z,
s,t)$ is a bounded, measurable function in all of
its variables for $|z|<\epsilon, |s|<\epsilon,
|t|<\epsilon$ (for details, see
e.g. \cite{BERbook},
\S 1.7.).

The Edge--of--the--Wedge theorem
ensures that this will be the case if $h$
extends holomorphically in $s+it$ for each
fixed $z$ to
$W_\Gamma^\epsilon$ and
$W_{-\Gamma}^\epsilon$ (for details, see
Section 2 of
\cite{BJT}) in a way such that $h$ is a
continuous function of all its variables in
$W_\Gamma^\epsilon\cup W_{-\Gamma}^\epsilon\cup
(M\cap \{|z|<\epsilon, |s|<\epsilon\})$. Actually,
this assumption might seem unnecessarily strong to
the reader; in fact, it is enough to assume that the
extension has boundary value in the sense of
distributions
$h$, but a regularity theorem implies that if
$h$ is continuous on $M$, then this
extension is actually continuous in the
sense above (see e.g. \cite{BERbook}, \S 7.3. and
\cite{Rosay}). Since we shall only deal with
functions which are at least continuous on
$M$, we shall assume continuity of the
extension a priori, and shall not go into
further detail.

For
the following arguments, we will fix
$\Gamma$ and say that
$h$ {\it extends up} differentiably (respectively
continuously)
if it extends holomorphically to
$W_\Gamma^\epsilon$ in a way such that the
extended function is a $C^1$ (respectively
continuous) function of all of its
variables in $\bar
W_\Gamma^\epsilon=W_\Gamma^\epsilon \cup M$
(where $M$ is shrunk appropriately) and that
$h$ {\it extends down} if it extends
holomorphically to
$\bar W_{-\Gamma}^\epsilon$ in a way such that the
extended function is a $C^1$ (respectively
continuous) function of all of its
variables (this is the
terminology  used in \cite{BERbook}, \S 9.2.,
in the case of hypersurfaces).

The main point can be phrased as follows:
\begin{center}
{\em If $h$ extends up continuously
then
$\bar h$ extends down continuously.}
\end{center}
In fact, the extension of $\bar{h}$ is just
given by $\bar h (z,\bar z, s, -t) =
\overline{h(z,\bar z, s, t)}$. This 
will be used later on for the components of
the CR--mapping $H$. We also need to know how
$Lh$ behaves for vector fields $L$
tangent to $M$ whose coefficients extend up
continuously. Here the result is that if the
coefficients of $L$ extend up continuously
and $h$ extends up differentiably, then $Lh$
extends up continuously. Similarily, if the
coefficients of $L$ extend up smoothly, the
regularity of the extension is dropped by $1$
(that is, if $h$ extends up in a $C^k$
manner, $Lh$ does so in a $C^{k-1}$ manner).
We summarize the discussion:
\begin{lem}\label{L:prel}  Let
$h$ be a CR--function on
$M$. If $h$ extends up continuously
(respectively differentiably), then $\bar h$
extends down continuously (respectively
differentiably). If $L$ is a vector field
tangent to $M$ whose coefficients extend up
smoothly and $h$ extends up of order $C^k$
then
$Lh$ extends up of order $C^{k-1}$. If $h$
extends up and down continuously, then $h$
extends holomorphically to a full
neighbourhood of the origin.
\end{lem}

The next lemma is used in order to actually
calculate with Definition~\ref{D:nondeg}. Its
proof is an easy induction on $k$ which is
left to the reader.
\begin{lem} \label{L:ekinbasis}
Let $L_1,\dots
, L_n$ be a local basis near $0$ of
$\Gamma(M,\crb{M})$.  For every multiindex
$\alpha=(\alpha_1,\dots,\alpha_n)\in\Nn$
define
${L}^\alpha =
{L_1}^{\alpha_1}
\cdots{L_n}^{\alpha_n}$.
Then
\begin{equation} \label{E:ekinbasis}
E_k(0) = \spanc
\{ {L}^\alpha
\rho^\prime_{l Z^\prime}(H(Z),\overline{H(Z)}
)
\big|_{Z=0}
\colon 1\leq l\leq d^\prime , |\alpha |\leq k \},
\end{equation}
where the $E_k(0)$ are defined by \eqref{E:ek}.
\end{lem}
\section{Proof of the Reflection Principle}
Choose coordinates
$(Z_1^\prime,\dots,Z_{N^\prime}^\prime)\in\CNp$
and a real--analytic defining function
$\rho^\prime$ for $M^\prime$ defined in some
neighbourhood of $0$. In these coordinates,
$H=(H_1,\dots,H_{N^\prime})$ with each $H_j$ a
CR--function on $M$. We will think of $M$ as
$\Cn\times\Rd$ as explained in
Section~\ref{S:prel}. A local basis of the
CR--vector fields (for which we will use the vector
notation and write 
$\Lambda=(\Lambda_1,\dots,\Lambda_n$)) is  then
given in matrix notation by  (with
$\varphi$ as in
\eqref{E:definingequM}; we also refer the
reader to \cite{BERbook}, \S 1.6.)
\begin{equation} \label{E:thelambdas}
\Lambda = \vardop{}{\bar z} - i
(\varphi_{\bar z})^T (I + i \varphi_s)^{-1}
\vardop{}{s}.
\end{equation}
The coefficients of all of these vector fields clearly
extend up and down (smoothly), since $\varphi$ is
real--analytic. Since
$H$ maps
$M$ into $M^\prime$, $\rho^\prime
(H(Z),\overline{H(Z)}) = 0$ for $Z\in M$.
Applying
$\Lambda_1,\dots,\Lambda_n$ repeatedly and
using the chain rule we see that for every
multiindex $\alpha\in\Nn$, $|\alpha|\leq
k_0$ (as in Lemma~\ref{L:ekinbasis},
$\Lambda^\alpha=
\Lambda^{\alpha_1}\cdots\Lambda^{\alpha_n}$)
and
for every
$l$,
$1\leq l \leq d^\prime$,
\begin{equation}\label{E:phialpha}
0 = \Lambda^\alpha \rho^\prime_l
(H(Z),\overline{H(Z)}) = \Phi_{l \alpha}
(H(Z),\overline{H(Z)}, (\Lambda^\beta
\overline{H(Z)})_{1\leq |\beta |\leq k_0} ), \quad
Z\in M,
\end{equation}
where $\Phi_{l\alpha}$ is a real--analytic
function defined and convergent on a
neighbourhood of $\{0\}\times\{0\}\times
\C^{K(k_0)}\subset\CNp\times
\CNp\times\C^{K(k_0)}$, $K(k_0)$ being the
cardinality of the set $\{ \beta\in\Nn\colon
1\leq |\beta | \leq k_0\}$. By our assumption and
Lemma~\ref{L:ekinbasis} we can choose
$N^\prime$ multiindices
$\alpha^1,\dots,\alpha^{N^\prime}$,
$0\leq |\alpha^j |\leq k_0$ and integers
$l^1,\dots,l^{N^\prime}$, $1\leq l^j \leq
d^\prime$, such that
\begin{equation}\label{E:fullspan} \spanc \{
\Lambda^{\alpha^j}
\rho^\prime_{l^j Z^\prime} (H(Z),\overline{H(Z)})
\big|_{Z=0} \colon 1\leq j \leq N^\prime
\} = \CNp.
\end{equation}
We consider the system of
equations
\begin{equation}\label{E:system}
\Phi_{l^j \alpha^j}
(X_1,\dots,X_{N^\prime},Y_1,
\dots,Y_{N^\prime},W)=0,\quad 1\leq j\leq
N^\prime,
\end{equation}
where $W\in\C^{K(k_0)}$. We claim that
\eqref{E:system} admits a (unique)
real--analytic solution in
$(X_1,\dots,X_{N^\prime})$ in a
neighbourhood of the point $(0,0,(\Lambda^\beta
\bar H (0))_{1\leq |\beta |\leq k_0})$. In fact, if
we compute the Jacobian of this system with
respect to $X_1,\dots, X_{N^\prime}$ it is of
full rank at this point because of
\eqref{E:fullspan}. So we can invoke the
implicit function theorem to conclude that
there are real--analytic functions
$\Upsilon_1,\dots,\Upsilon_{N^\prime}$,
convergent on a neighbourhood $U$ of
$(0,(\Lambda^\beta
\bar H (0))_{1\leq |\beta |\leq k_0})$, such that
the unique solution of \eqref{E:system} in
$U$ is given by
\begin{equation}\label{E:sol}
X_j = \Upsilon_j (Y_1,\dots,Y_{N^\prime},W),
\quad 1\leq j \leq N^\prime.
\end{equation}
Recalling \eqref{E:phialpha} we conclude that
\begin{equation}\label{E:solinH}
H_j (Z) = \Upsilon_j (\overline{H_1 (Z)},
\dots, \overline{H_{N^\prime}
(Z)},(\Lambda^\beta \overline{H(Z)} )_{1\leq |\beta
|\leq k_0 } ).
\end{equation}
Now the proof is finished by using
Lemma~\ref{L:prel}: Each $H_j$ is assumed to
extend up; so each $\bar{H}_j$ extends down
(of order $C^{k_0}$, by a regularity theorem,
see Theorem 7.5.1. in \cite{BERbook}). Hence
the whole right hand side of
\eqref{E:solinH}  extends down
continuously (after choosing $\epsilon$ small
enough). This shows that each $H_j$ extends
down continuously. Since each $H_j$ also extends up,
Lemma~\ref{L:prel} implies that each $H_j$
extends holomorphically to a full
neighbourhood of the origin. The theorem is
proved.
\bibliographystyle{plain}
\bibliography{refl2}
\end{document}